\documentclass[12pt]{article}
\usepackage{amsfonts,amsmath,amsxtra}
\usepackage{latexsym}
\usepackage[matrix,arrow,curve]{xy}
\usepackage{amssymb}
\def\hybrid{\topmargin 0pt      \oddsidemargin 0pt
        \headheight 0pt \headsep 0pt
        \textwidth 16.5cm
        \textheight 23cm
        \marginparwidth 0.0in
        \parskip 5pt plus 1pt   \jot = 1.5ex}
\catcode`\@=11
\def\marginnote#1{}
\newcount\hour
\newcount\minute
\newtoks\amorpm
\hour=\time\divide\hour by60 \minute=\time{\multiply\hour by60
\global\advance\minute by-\hour}
\edef\standardtime{{\ifnum\hour<12 \global\amorpm={am}%
        \else\global\amorpm={pm}\advance\hour by-12 \fi
        \ifnum\hour=0 \hour=12 \fi
      \number\hour:\ifnum\minute<10 0\fi\number\minute\the\amorpm}}
\edef\militarytime{\number\hour:\ifnum\minute<10 0\fi\number\minute}

\def\draftlabel#1{{\@bsphack\if@filesw {\let\thepage\relax
   \xdef\@gtempa{\write\@auxout{\string
      \newlabel{#1}{{\@currentlabel}{\thepage}}}}}\@gtempa
   \if@nobreak \ifvmode\nobreak\fi\fi\fi\@esphack}
        \gdef\@eqnlabel{#1}}
\def\@eqnlabel{}
\def\@vacuum{}
\def\draftmarginnote#1{\marginpar{\raggedright\scriptsize\tt#1}}

\def\draft{\oddsidemargin -0.1truein
        \def\@oddfoot{\sl preliminary draft \hfil
        \rm\thepage\hfil\sl\today\quad\militarytime}
        \let\@evenfoot\@oddfoot \overfullrule 3pt
        \let\label=\draftlabel
        \let\marginnote=\draftmarginnote
\def\@eqnnum{{\rm (\theequation)}
\rlap{\kern\marginparsep\tt\@eqnlabel}%
\global\let\@eqnlabel\@vacuum}  }

\newfont{\Bbbb}{msbm7 scaled 1\@ptsize00}
\newcommand{\zs}{\raise-1pt\hbox{$\mbox{\Bbbb Z}$}}

\def\numberbysection{\@addtoreset{equation}{section}
        \def\theequation{\thesection.\arabic{equation}}}
\numberbysection

\renewcommand{\theequation}{\thesection.\arabic{equation}}

\def\titlepage{\@restonecolfalse\if@twocolumn\@restonecoltrue\onecolumn
     \else \newpage \fi \thispagestyle{empty}\c@page\z@
\def\thefootnote{\fnsymbol{footnote}} }
\def\endtitlepage{\if@restonecol\twocolumn \else  \fi
        \def\thefootnote{\arabic{footnote}}
        \setcounter{footnote}{0}}  
\relax
\hybrid
\makeatletter
\newdimen\normalarrayskip            
\newdimen\minarrayskip               
\normalarrayskip\baselineskip \minarrayskip\jot
\newif\ifold             \oldtrue            \def\new{\oldfalse}
\def\arraymode{\ifold\relax\else\displaystyle\fi}
\def\eqnumphantom{\phantom{(\theequation)}} 
\def\@arrayskip{\ifold\baselineskip\z@\lineskip\z@
     \else
     \baselineskip\minarrayskip\lineskip1\baselineskip\fi}
\def\@arrayclassz{\ifcase \@lastchclass \@acolampacol \or
\@ampacol \or \or \or \@addamp \or
   \@acolampacol \or \@firstampfalse \@acol \fi
\edef\@preamble{\@preamble
  \ifcase \@chnum
     \hfil$\relax\arraymode\@sharp$\hfil
     \or $\relax\arraymode\@sharp$\hfil
     \or \hfil$\relax\arraymode\@sharp$\fi}}
\def\@array[#1]#2{\setbox\@arstrutbox=\hbox{\vrule
     height\arraystretch \ht\strutbox
     depth\arraystretch \dp\strutbox
width\z@}\@mkpream{#2}\edef\@preamble{\halign \noexpand\@halignto
\bgroup \tabskip\z@ \@arstrut \@preamble \tabskip\z@ \cr}%
\let\@startpbox\@@startpbox \let\@endpbox\@@endpbox
  \if #1t\vtop \else \if#1b\vbox \else \vcenter \fi\fi
  \bgroup \let\par\relax
  \let\@sharp##\let\protect\relax
  \@arrayskip\@preamble}
\def\eqnarray{\stepcounter{equation}%
              \let\@currentlabel=\theequation
              \global\@eqnswtrue
              \global\@eqcnt\z@
              \tabskip\@centering              
              \let\\=\@eqncr
              $$%
            \halign to \displaywidth  \bgroup
             \eqnumphantom \@eqnsel
      \hskip\@centering                               
    $\displaystyle  \tabskip\z@ {##}$%
    &\global\@eqcnt\@ne \hskip 2\arraycolsep
         $ \displaystyle  \arraymode{##}$\hfil
    &\global\@eqcnt\tw@ \hskip 2\arraycolsep
         $\displaystyle\tabskip\z@{##}$\hfil
         \tabskip\@centering
    &{##}\tabskip\z@\cr}
\makeatother

\def\IC{\mathbb{C}}

\def\IF{\mathbb{F}}

\def\IR{\mathbb{R}}
\def\IZ{\mathbb{Z}}

\def\CW {\mathcal{W}}

\def\Fg{{\frak g}}

\def\fg{\mathfrak{g}}
\def\fh{\mathfrak{h}}

\def\a {{\alpha}}

\def\g {{\gamma}}
\def\s {{\sigma}}

\def\vs{\varsigma}
\def\vsb{\bar{\varsigma}}

\def\Ad{{\mathop{\rm Ad}}}

\def\Lie{{\rm Lie}}

\def\frak{\mathfrak}

\def\<{\langle}
\def\>{\rangle}
\def\s{\sigma}
\def\sb{\bar{\sigma}}

\newtheorem{te}{Theorem}[section]
\newtheorem{de}{Definition}[section]
\newtheorem{prop}{Proposition}[section]           

\newtheorem{lem}{Lemma}[section]
\newtheorem{ex}{Example}[section]

\newcommand\bqa{\begin{eqnarray}}
\newcommand\eqa{\end{eqnarray}}
\def\be{\begin{eqnarray}\new\begin{array}{cc}}
\def\ee{\end{array}\end{eqnarray}}
\def\beq{\begin{equation}}
\def\eeq{\end{equation}}
\def\bse{\begin{subequations}}                
\def\ese{\end{subequations}}
\def\bp{\begin{pmatrix}}
\def\ep{\end{pmatrix}}

\newcommand{\proof}{\noindent {\it Proof}. }
\newcounter{pac}[section]

\newcounter{pacc}[subsection]

0

\begin{document}


\title{\bf Normalizers of maximal tori \\ and  real forms of Lie groups}
\author{A.A. Gerasimov, D.R. Lebedev and S.V. Oblezin}
\date{}
\maketitle

\renewcommand{\abstractname}{}

\begin{abstract}
\noindent {\bf Abstract}.  Given a complex connected reductive Lie
group $G$ with a maximal torus  $H\subset G$,  Tits defined an
extension $W_G^T$ of the corresponding Weyl group $W_G$. The
extended group is supplied with an embedding into the normalizer
$N_G(H)$, such that $W_G^T$ together with $H$ generate $N_G(H)$. In
this paper  we propose an interpretation of the Tits classical
construction in terms of the maximal split real form $G(\IR)\subset
G$, which leads to the simple topological description of $W^T_G$. We
also consider a variation of the Tits construction associated with
compact real form $U$ of $G$. In this case we define an extension
$W_G^U$ of the Weyl group $W_G$, naturally embedded into the group
extension $\widetilde{U}:=U\rtimes\Gamma$ of the  compact real form
$U$ by the Galois group $\Gamma={\rm
  Gal}(\IC/\IR)$. Generators of $W^U_G$ are squared to identity as in the Weyl
group $W_G$.  However, the non-trivial action of $\Gamma$ by outer
automorphisms requires $W^U_G$ to be a non-trivial extension of
$W_G$. This gives a specific presentation of the maximal torus
normalizer of the group extension $\widetilde{U}$. Finally, we
describe explicitly the adjoint action of $W_G^T$ and $W^U_G$ on the
Lie algebra of $G$.

\end{abstract}


\vspace{5 mm}

\section{Introduction}

In the standard approach to classification of complex semisimple Lie
groups the problem is reduced to an  equivalent problem of
classification of root data. In other words, the root data (the
roots/co-roots and the lattices of characters/co-characters of a
maximal torus $H\subset G$) define the corresponding semisimple Lie
group up to an isomorphism. Curtis, Wiederhold and Williams
\cite{CWW} demonstrate that for classification of compact connected
semisimple Lie groups $G$ it is enough to classify the normalizers
$N_G(H)$ of  maximal tori $H\subset G$. The normalizer provides an
information about the action of the Weyl group $W_G:=N_G(H)/H$ on
$H$, but this is not enough for classification, as one also needs
the precise structure of the group extension of $W_G$ by $H$. Thus,
for the classification problem one might replace the original
object, the semisimple Lie group $G$, by the group extension of the
finite group $W_G$ by the abelian Lie group $H$. One perspective to
grasp this equivalence is to look at $N_G(H)$ as a kind of
degeneration of $G$ \cite{CWW}. An apparently related but more
conceptual approach is based on attempts to look at $N_G(H)$ as the
Lie group $G$ defined over some non-standard  base field
 (akin to the mysterious field $\IF_1$ ``with one element'' introduced
by Tits \cite{T3}, probably with regard to this subject). In this
way an equivalence of the two classification problems for compact
semisimple Lie groups and normalizers looks like a manifestation of
a general principle (due to C. Chevalley \cite{C}), saying that a
classification of semisimple algebraic groups should  not
essentially depend  on the nature of the base algebraically closed
field.

The reasoning above demands a more detailed study of the group
extension structure on $N_G(H)$. The key fact is that this extension
does not split in general \cite{D}, \cite{T1}, \cite{T2},
\cite{CWW}, \cite{AH}. To get a universal description of $N_G(H)$
one should look for a set-theoretic section of the projection
$N_G(H)\to W_G$  generating an  extension of $W_G$. Such
construction was proposed by Demazure \cite{D}  and Tits \cite{T1},
\cite{T2}. It may be naturally formulated in terms of the Tits
extension $W^T_G$ of Weyl group  $W_G$ by the subgroup of order $2$
elements in $H$.  This construction allows an explicit presentation
of $N_G(H)$ by generators and relations.

Although the Tits construction is known for a long time,
  its forthright explanation involves scheme-theoretic arguments. Precisely,
for a Chevalley group scheme $\mathcal{G}$ over $\IZ$ the Tits
extension is the group of $\IZ$-points of the normalizer of the
$\IZ$-split torus of $\mathcal{G}$. In this paper we use
set-theoretic arguments to explain the Tits construction in the case
of complex reductive Lie groups (for recent discussion on Tits
groups see e.g. \cite{N}, \cite{DW}, \cite{AH}). After reminding
general results on  normalizers of maximal tori in Section 2 we
revisit the Tits construction in Section 3. We stress that the Tits
group construction is defined for the split real form $G(\IR)\subset
G$ of a complex semisimple group $G$. This enables us to present in
Proposition \ref{sdot}  a simple purely topological description of
the Tits extension of the Weyl group $W_G$ (our considerations
appear to be pretty close to the final section of \cite{BT}).

As a variation of the Tits
construction we consider its analog for the
maximal compact group $U\subset G$. Precisely, we define an extension
$W_G^U$ of $W_G$ embedded into a semi-direct product $\widetilde{U}=U\rtimes\Gamma$
of the maximal compact subgroup $U\subset G$ and the
 Galois group $\Gamma={\rm Gal}(\IC/\IR)$ acting in the standard
(non-twisted) way on $U$. In this case the natural generators of the
corresponding extension $W_G^U$ of the Weyl group are squared to
identity (in contrast with the case of the Tits group), while the
non-trivial extension of $W_G$ arises from the action of $\Gamma$
via on $W_G^U$ by outer automorphisms. The main result of this paper
(Theorem \ref{TE} in Section 4) describes the structure of maximal
tori normalizers in the Galois group extension $\widetilde{U}$ of
the compact connected semisimple Lie group $U$. Note that our
construction is dealing with the extended group $\widetilde{U}$ and
thus, it differs from the constructions given in \cite{DW}, \cite
{AH}. In Section 5 we calculate explicitly the adjoint action of the
Tits group and of its unitary analog
 on the Lie algebra  $\Fg={\rm
  Lie}(G)$. This action, in contrast with the adjoint
action on $\mathfrak{h}\subset \Fg$,  depends on the lift of $W_G$
into $G$. Finally in Section 6 we provide details of the proof of Theorem \ref{TE}.

{\it  Acknowledgements:} We are grateful to the referees for careful
reading and valuable suggestions including a short proof of Lemma
\ref{Dihedral}. The research of the second author was supported by
the RSF grant 16-11-10075. The work of the third author was
partially supported by the EPSRC grant EP/L000865/1 and by the RSF
grant 16-11-10075. The third author is also thankful to the Max
Planck Institute for Mathematics in Bonn, where his work on the
project was started.

\section{Normalizers of maximal tori and Weyl groups}

We start with recalling standard facts on  normalizers of maximal
tori and the associated Weyl groups. Let $G$ be a complex connected
semisimple Lie group, $H\subset G$ be a maximal torus and $N_G(H)$
be its normalizer in $G$. Then there is the following exact sequence
 \be\label{NGH}
 1\longrightarrow H\longrightarrow N_G(H)\stackrel{p}\longrightarrow
 W_G\longrightarrow 1,
 \ee
where $p$ is the projection on the finite group $W_G:=N_G(H)/H$, the
Weyl group of $G$. The Weyl group $W_G$ does not actually depend on
the choice of $H\subset G$, and thus appears to be an invariant of
$G$. Let $\Fg:={\rm Lie}(G)$ and let $I$ be the set of vertices of
the Dynkin diagram associated to $G$, where $|I|={\rm rank}(\Fg)$.
Let $(\Delta,\Delta^\vee)$ be the root-coroot system corresponding
to $G$, $\{\alpha_i,\,i\in I\}$ be a set of positive simple roots,
and $\{\alpha_i^\vee,\,i\in I\}$ be the corresponding set of
positive simple co-roots. Let $\|a_{ij}\|$,
$a_{ij}=\<\alpha_i^\vee,\alpha_j\>$, be the Cartan matrix of
$(\Delta,\,\Delta^\vee)$. The Weyl group $W_G$ has a simple
presentation in terms of generators and relations. Precisely, $W_G$
is generated by simple root reflections $\{s_i,\,i\in I\}$ subjected
to
 \be\label{BR0}
  s_i^2=1,
 \ee
 \be\label{BR}
 \underbrace{s_is_js_i\cdots }_{m_{ij}}\,
 =\,\underbrace{s_js_is_j\cdots }_{m_{ij}}\,,\qquad i\neq j\,,\quad i,j\in I\,,
 \ee
where $m_{ij}=2,3,4,6$ for $a_{ij}a_{ji}=0,1,2,3$, respectively.
Equivalently these relations   may be written in the Coxeter form:
 \be
  s_i^2=1, \qquad  (s_is_j)^{m_{ij}}=1\,, \qquad i\neq j\,,\quad i,j\in I\,.
 \ee
The exact sequence \eqref{NGH} defines the canonical action of $W_G$
on $H$. Let $h_i\in\mathfrak{h}={\rm Lie}(H)$ be the generators
corresponding to the simple co-roots $\a_i^{\vee}$, then the
$W_G$-action on $\fh\subset\fg$ and on its dual is as follows
 \be\label{WeylGroupAction}
  s_i(h_j)=h_j-\<\a_i,\a_j^\vee\>h_i\,=\,h_j-a_{ji}h_i\,,\\
 s_i(\a_j)=\a_j-\<\a_j,\a_i^\vee\>\a_i\,=\,\a_j-a_{ij}\a_i\,.
 \ee
The exact sequence  \eqref{NGH} does not split in
 general, i.e. $N_G(H)$ is not necessarily isomorphic to a semi-direct product of $W_G$ and
 $H$. A delicate situation in this regard is described by the following
result due to  \cite{CWW}, \cite{AH}.
\begin{te} \label{BRTH}  Assume $G$ is a simple complex Lie group and let $Z(G)$ be
 the  center of $G$. Then the
exact sequence \eqref{NGH} splits in the following and only the
following cases:
\begin{itemize}
\item Type $A_{\ell},\,\ell\geq1,\,\ell\neq3$, such that  $|Z(G)|$ is odd;

\item Type $B_{\ell},\,\ell>1$, for the adjoint form;


\item Type $D_{\ell},\,\ell>2$, for all forms except $Spin(2\ell)$;

\item Type $G_2$.
\end{itemize}
\end{te}

Thus to have an explicit description of the normalizer $N_G(H)$ one
should  look for appropriate   set-theoretic section of the
projection map $p$ in \eqref{NGH} generating some extension of
  $W_G$. In the following  section we
provide the construction of the resulting extension of the Weyl
group by a finite group. Let us note that for a normal finite
subgroup $Z\subset G$ one has: if \eqref{NGH} splits for $G$ then it
splits for $G/Z$ as well. In the following,   for simplicity,  we
consider only the case of simply-connected complex groups.

\section{The Tits extension of Weyl group}

To describe the extension \eqref{NGH} in terms of generators and
relations Tits proposed the following
extension $W^T_G$ of the Weyl group $W_G$ by a discrete group \cite{T1},
\cite{T2} (closely related results were obtained by Demazure \cite{D}).

\begin{de} Let $A=\|a_{ij}\|$ be the Cartan matrix corresponding to a
  semisimple Lie algebra $\Fg={\rm Lie}(G)$  and let
$m_{ij}=2,3,4,6$ for $a_{ij}a_{ji}=0,1,2,3$, respectively. The Tits
group $W^T_G$ is  an extension of the Weyl group
$W_G$ by an abelian group $\IZ_2^{|I|}$
generated by $\{\tau_i,\,\theta_i,\,i\in I\}$
subjected to the following relations:
  \be\label{TitsSq}
 (\tau_i)^2\,=\,\theta_i\,,\qquad
 \theta_i\theta_j=\theta_j\theta_i\,, \qquad\theta_i^2\,=\,1\,,
 \ee
 \be\label{TitsAd}
  \tau_i\theta_j\,=\,\theta_i^{-a_{ji}}\theta_j\tau_i\,,
 \ee
 \be\label{Tits1}
  \underbrace{\tau_i\tau_j\tau_i\cdots }_{m_{ij}}\,
  =\,\underbrace{\tau_j\tau_i\tau_j\cdots }_{m_{ij}}\,\,,\qquad i\neq j\,,\quad i,j\in I\,,
  \ee
where the abelian subgroup is generated by $\{\theta_i,\,i\in I\}$.
\end{de}

Let $\{h_i,\,e_i,\,f_i\,: i\in I\}$ be the Chevalley-Serre
generators of the  Lie algebra $\Fg={\rm Lie}(G)$, satisfying the
standard relations
 \be\label{Serre1}
  [h_i,\,e_j]\,=\,a_{ij}e_j\,,\qquad[h_i,\,f_j]\,=\,-a_{ij}f_j\,,\qquad
  [e_i,\,f_j]\,=\,\delta_{ij}h_j\,,
  \ee
  \be\label{Serre2}
  {\rm ad}_{e_i}^{1-a_{ij}}(e_j)=0, \qquad   {\rm
    ad}_{f_i}^{1-a_{ij}}(f_j)=0,
 \ee
where $A=\|a_{ij}\|$ is the Cartan matrix i.e.
$a_{ij}=\<\a_i^\vee,\,\a_j\>$.

 According to \cite{BT} (see also \cite{T1}) there exists a subset
$\{\zeta_i,\,i\in I\}\subset H$ of canonical elements of order two
satisfying the following relations
 \be
  s_i(\zeta_j)=\zeta_j\,\zeta_i^{-a_{ji}}\,,
  \ee
where $s_i$, $i\in I$  are the generators of the Weyl group $W_G$
\eqref{BR0},  \eqref{BR}.

 \begin{te}[Demazure-Tits] \label{DTT}
   Let $W^T_G$ be the Tits group associated with the complex
semisimple simply connected Lie group $G$, then the map
 \be\label{Tits}
  \tau_i\longmapsto\dot{s}_i:=e^{f_i}\,e^{-e_i}\,e^{f_i}, \qquad
  \theta_i\longmapsto\zeta_i\,,\qquad i\in I\,,
 \ee
defines an embedding of the Tits group $W^T_G$ into $N_G(H)$, such
that $p(W_G^T)=W_G$ for the projection $p$ in \eqref{NGH}. In
particular, the normalizer group $N_G(H)$ is generated by $H$ and by
the image of the Tits group, so that the following relations hold:
 \be
 \Ad_{\dot{s}_i}(h)\,=\,s_i(h)\,,\qquad\forall
 h\in\mathfrak{h}=\Lie(H)\,,  \qquad i\in I.
 \ee
\end{te}

\begin{ex}  In the standard faithful two-dimensional representation $\phi:
  SL_2(\IC)\to {\rm End}(\IC^2)$  given by \eqref{STREP} we have
  \be
  \phi(\dot{s})\,=\,
  \begin{pmatrix} 0 & -1 \\ 1 & 0\end{pmatrix},\qquad
    \phi(\dot{s})^2=\phi(\zeta)=\begin{pmatrix} -1 & 0 \\ 0 & -1\end{pmatrix}.
\ee
\end{ex}

The appearance of the Tits extension $W_G^T$ via a specific choice
of a set-theoretic section  looks  a bit  ad hoc.  As it is
mentioned in Introduction
  one may use a scheme-theoretic argument to support this particular
  choice of extension of $W_G$.
  In the following we propose a set-theoretic
  argument based on consideration of
  the split real form $G(\IR)$ of $G$
 to elucidate the construction of $W_G^T$.
For the split real form $G(\IR)\subset G$ there is an analog of
\eqref{NGH}:
 \be\label{NGHR}
  1\longrightarrow
  H(\IR)\longrightarrow N_{G(\IR)}(H(\IR)) \stackrel{p}\longrightarrow
  W_G\longrightarrow 1,
 \ee
with the real split maximal torus  given by the intersection
 \be
  H(\IR)=H\cap G(\IR),
 \ee
of the complex maximal torus  with the split real subgroup. Note
that the set-theoretic section of \eqref{NGHR} defining $W_G^T$
provides a set-theoretic section of \eqref{NGH} thus embedding
$W_G^T$ into $N_{G(\IR)}(H(\IR))$.  The group $H(\IR)$ allows the
product decomposition
 \be
  H(\IR)=MA, \qquad M:=H(\IR)\cap K,
 \ee
where $K\subset G(\IR)$ is a maximal compact subgroup of $G(\IR)$,
$M$ is isomorphic to the group $\IZ_2^{|I|}$ and $A\subset H(\IR)$
is the connected exponential group $A=\exp(\fh(\IR))$ without
torsion. Therefore, $H(\IR)$ is not connected and consists of
$2^{|I|}$ components, and the group $M$ may be identified with the
discrete group of connected components of $H(\IR)$:
 \be
  M\,=\,\pi_0(H(\IR))\,.
 \ee
Considering the groups of connected components of the topological
groups entering \eqref{NGHR} we obtain  the induced exact sequence
 \be\label{NGH1}
  1\longrightarrow \pi_0(H(\IR))\longrightarrow \pi_0(N_{G(\IR)}(H(\IR)))
  \longrightarrow W_G\longrightarrow 1\,,
 \ee
so that $|\pi_0(N_{G(\IR)}(H(\IR)))|=|M|\cdot|W_G|=|W_G^T|$. This
provides a canonical extension of $W_G$ by
  the abelian group $M$ of order $2^{|I|}$.

 Explicitly, the group  of connected components may be identified
with the quotients by the connected normal subgroup $A$
 \be\label{NGHR1}
  \pi_0(N_{G(\IR)}(H(\IR))) \simeq N_{G(\IR)}(H(\IR)/A\,,
 \ee
and we have the exact sequence
 \be\label{EXT0}
  1\longrightarrow A\longrightarrow N_{G(\IR)}(H(\IR)) \longrightarrow
  \pi_0(N_{G(\IR)}(H(\IR)))\longrightarrow 1\,.
 \ee

\begin{lem}\label{SPLIT} The exact sequence \eqref{EXT0} splits and thus
$\pi_0(N_{G(\IR)}(H(\IR)))$ allows an embedding into
$N_{G(\IR)}(H(\IR))$.
  \end{lem}

  \proof The extension \eqref{EXT0} is an instance of extensions of
$\pi_0(N_{G(\IR)}(H(\IR)))$ by $A$. Such extensions  are classified
  by the group $H^2(\pi_0(N_{G(\IR)}(H(\IR))),A)$. The
  triviality of this group follows from the fact that
$\pi_0(N_{G(\IR)}(H(\IR)))$ is a finite group and   $A$ is an
instance of abelian torsion-free group. By the standard cohomology
argument (see e.g. \cite{B}, Chapter VI), the second
  cohomology  of any finite group
with coefficients in a torsion-free abelian group is trivial. Thus
the extension \eqref{EXT0} is necessarily trivial and therefore the
required embedding exists. $\Box$

Up to now we have constructed a canonical extension of $W_G$ given
by \eqref{NGH1}. It is easy to see that this extension is isomorphic
to
  the Tits group.

\begin{prop}\label{sdot}
The following isomorphism holds
 \be \label{pi0}
  \pi_0(N_{G(\IR)}(H(\IR)))\,\simeq\, W^T_G\,.
 \ee
\end{prop}

\proof Recall that  the images $\dot{s}_i$, $\zeta_i$, $i\in I$ of
the Tits generators belong to  the maximally split real subgroup
$G(\IR)\subset G$. Thus we have the embedding $W_G^T$ into
$N_{G(\IR)}(H(\IR))$ providing section of the exact sequence
 \be\label{NWText}
 1\longrightarrow A\longrightarrow N_{G(\IR)}(H(\IR))
 \longrightarrow W^T_G\longrightarrow 1\,.
 \ee
Note that splitting of this exact sequence may be independently
verified by using an argument similar to the one used in our proof
of Lemma \ref{SPLIT}. Therefore the group $N_{G(\IR)}(H(\IR))$
allows a representation as semidirect product
$N_{G(\IR)}(H(\IR))=A\rtimes W^T_G$. By considering the connected
components we deduce the assertion \eqref{pi0}. $\Box$

\begin{ex} For maximal split form $SL_2(\IR)\subset SL_2(\IC)$  we have
  \be
  H(\IR)=\left\{\begin{pmatrix} \lambda & 0\\ 0& \lambda^{-1}\end{pmatrix},\,
  \lambda\in \IR^*\right\}, \qquad H(\IR)=MA\,,
  \ee
  \be
   A=\left\{\begin{pmatrix} \lambda & 0\\ 0& \lambda^{-1}\end{pmatrix},\,
  \lambda\in \IR^*_+\right\},\qquad
    M=\left\{\pm {\rm Id}\right\}\,.
    \ee
Elements $g\in N_{SL_2(\IR)}(H(\IR))$ are defined by the condition that for each
    $\lambda\in \IR^*$ there exists $\tilde{\lambda}\in \IR^*$ such that
 \be\label{NORM}
  g\begin{pmatrix}\lambda & 0
  \\ 0 & \lambda^{-1}\end{pmatrix} =\begin{pmatrix}\tilde{\lambda} & 0
  \\ 0 & \tilde{\lambda}^{-1}\end{pmatrix}g\,,
  \qquad g=\begin{pmatrix}a & b \\ c & d\end{pmatrix},\qquad ad-bc=1.
 \ee
One might check that the normalizer group $N_{SL_2(\IR)}(H(\IR))$ is
a union of two components:
 \be
  N_{SL_2}(H(\IR))=N_1\sqcup N_s\,,
 \ee
where $N_1$ is a set of diagonal elements with $c=b=0,\,ad=1\neq 0$,
and $N_s$ is the set of anti-diagonal elements with $a=d=0,\,cb=-1$.
Each of the co-sets $N_1,\,N_s$ splits further into two connected
components
 \be
  N_1=N^+_1\sqcup N^-_1, \qquad N_s=N^+_s\sqcup N^-_s,
 \ee
depending on the sign of the entries $c,d$ in the last row of $g$.

The group $\pi_0(N_{SL_2(\IR)}(H(\IR)))$ is isomorphic to the
quotient of $N_{SL_2(\IR)}(H(\IR))$ by  $A\simeq\IR^*_+$ and
consists of four  elements corresponding to the connected components
$N_1^{\pm}$, $N_s^{\pm}$ of the group $N_{SL_2}(H(\IR))$, allowing
the following parameterization:
 \be
  N_1^+=A,\qquad N_s^+=\dot{s}A, \qquad N_1^-=\theta A, \qquad N_s^-=\theta
  \dot{s} A\,,
 \ee
where
 \be
  \dot{s}=\begin{pmatrix} 0 & -1 \\ 1& 0\end{pmatrix}\,, \qquad
  \theta=(\dot{s})^2 =\begin{pmatrix} -1 & 0 \\0&- 1 \end{pmatrix}\,, \qquad
  \theta \dot{s} =\begin{pmatrix} 0 & 1\\ -1& 0 \end{pmatrix}\,.
 \ee
Clearly, the group $\pi_0(N_{SL_2(\IR)}(H(\IR)))$ is generated by
representatives of  the connected components $N_{1,s}^\pm$, so
indeed it is isomorphic to the order four cyclic group generated by
$\dot{s}$, in concordance with \eqref{pi0}.
\end{ex}

\section{Weyl group  and
  Galois extension of the compact real form}

As we have demonstrated in the previous Section  the Tits group
extension $W^T_G$ appears  quite naturally if we consider the split
real subgroup $G(\IR)\subset G$. This motivates to look for analogs
of the Tits construction associated with other real forms of $G$.

Let $\{h_i,\,e_i,\,f_i\,; i\in I\}$ be the Chevalley-Serre
generators of the  Lie algebra $\Fg={\rm Lie}(G)$, satisfying
\eqref{Serre1}, \eqref{Serre2}. We fix the split real structure
by assuming the generators $\{h_i,\,e_i,\,f_i\,; i\in I\}$ to be
real. Let $\top:\Fg\to\Fg$ be the Cartan anti-involution,
associated with the split real structure:
 \be\label{Cartaninv}
  \top\,:\quad e_i\,\longmapsto f_i\,,\quad
  f_i\,\longmapsto\,e_i\,,\quad h_i\,\longmapsto\,h_i\,,\qquad i\in
  I\,.
 \ee
Let  $U\subset G$ be the connected  compact  real form of the Lie
group $G$:
 \be
  U=\{g\in G\,:\,g^\dagger g=1\}\,,
 \ee
where $g\mapsto g^\dagger$ is the composition of the Cartan
anti-involution \eqref{Cartaninv} with the complex conjugation
associated with the split real structure).  The Galois group
$\Gamma:={\rm Gal}(\IC/\IR)\simeq\IZ_2$ of the extension $\IR\subset
\IC$ is generated by $\gamma$, $\gamma^2=1$. The group $\Gamma$ acts
both on $G$ and on $U$ by the complex conjugation, so let us
introduce the following semidirect products:
 \be\label{SDP}
  \widetilde{U}\,:=\,U\rtimes \Gamma\,\subset\,\widetilde{G}:=G\rtimes\Gamma\,,
 \ee
Since the generators $\{e_i,f_i,h_i\,;\,i\in I\}$ are real, they are
fixed by $\gamma\in\Gamma$. Let us  note that the $\Gamma$-fixed
subgroup of $U$ is a maximal compact subgroup $K\subset G(\IR)$ of
the split real form $G(\IR)$.

\begin{de}\label{DDD} Let $\|a_{ij}\|$ be the Cartan matrix corresponding to a
  semisimple Lie algebra $\Fg={\rm Lie}(G)$. Let $W^U_G$ be a group generated by
  $\{\sigma_i,\,\bar{\sigma}_i\,;\,i\in
I\}$ subjected to
 \be\label{rel}
  \s_i^2\,=\,\sb_i^2\,=\,1\,,\qquad
  \s_i\sb_i\,=\,\sb_i\s_i\,,
 \ee
 \be\label{rel0}
 \s_j\s_i\s_j\,=\,\sb_j\,\Pi^{1-a_{ji}}(\sb_i)\,\sb_j\,,\qquad i\neq j\,,\quad i,j\in I\,,
  \ee
  \be\label{Mrel}
  \underbrace{\s_i\s_j\s_i\cdots}_{m_{ij}}\,=\,
  \underbrace{\sb_j\sb_i\sb_j\cdots}_{m_{ij}}\,,
  \qquad i\neq j\,,\quad i,j\in I\,,
  \ee
where in \eqref{Mrel} $m_{ij}=2,3,4,6$ for $a_{ij}a_{ji}=0,1,2,3$, respectively.
 Here $\Pi$ is the involutive map of the generators given by $\Pi(\s_i)=\sb_i$,
$\Pi(\sb_i)=\s_i$, $i\in I$.
\end{de}
\begin{lem} For all $i,j\in I$ the following holds:
 \be\label{etacomm}
  (\s_j\sb_j)(\s_i\sb_i)=(\s_i\sb_i)(\s_j\sb_j)\,.
 \ee
\end{lem}
\proof For both $a_{ji}$ and $a_{ij}$ odd, \eqref{rel0} reads
 \be
  (\s_j\sb_j)\,\s_i(\s_j\sb_j)^{-1}=\sb_i\,,\qquad(\s_j\sb_j)\,\sb_i(\s_j\sb_j)^{-1}=\s_i,
 \ee
which yields $(\s_j\sb_j)(\s_i\sb_i)=(\s_i\sb_i)(\s_j\sb_j)$.

For $a_{ji}$ even, we have two cases: $a_{ji}=a_{ij}=0$ and
$a_{ji}=-2,\,a_{ij}=-1$. In the former case \eqref{Mrel} gives
$\s_i\s_j=\sb_j\sb_i$ and $\s_j\s_i=\sb_i\sb_j$, which implies
$$
 (\s_j\sb_j)(\s_i\sb_i)=\s_j(\sb_j\sb_i)\s_i=(\s_j\s_i)(\s_j\s_i)
 =\sb_i(\sb_j\sb_i)\sb_j=(\sb_i\s_i)(\s_j\sb_j)\,.
$$
In the case $a_{ji}=-2,\,a_{ij}=-1$ \eqref{rel0} gives
$(\s_j\sb_j)\s_i=\s_i(\s_j\sb_j)$ and
$(\s_i\sb_i)\s_j=\sb_j(\s_i\sb_i)$, so that the latter equality
entails $\sb_j(\s_i\sb_i)=(\s_i\sb_i)\s_j$. Thus we have
$$
 (\s_i\sb_i)(\s_j\sb_j)=\sb_j(\s_i\sb_i)\sb_j=(\s_j\sb_j)(\s_i\sb_i)\,.
$$
This completes our proof. $\Box$

\begin{lem}\label{INV} The map $\Pi$ acting on the generators by
  \be\label{inviso}
  \Pi(\s_i)=\sb_i, \qquad \Pi(\sb_i)=\s_i\,,
  \ee
extends to an involutive  automorphism of the group $W^U_G$.
\end{lem}

\proof   Clearly, the
relations \eqref{rel}, \eqref{Mrel} are invariant under the action of
$\Pi$. The relation \eqref{rel0} transforms under the $\Pi$-action
into
 \be\label{rel0bar}
 \sb_j\sb_i\sb_j\,=\,\s_j\,\Pi^{1-a_{ji}}(\s_i)\,\s_j\,,\qquad i\neq
 j\,,  \quad i,j\in I\,,
 \ee
which may  be equivalently written as follows
 \be\label{rel0bar1}
  (\s_j\sb_j)\,\sb_i(\s_j\sb_j)^{-1}\,=\,\Pi^{-a_{ji}}(\sb_i)\,,\qquad
  i\neq j\,, \quad i,j\in I\,.
  \ee

For $a_{ji}$ odd \eqref{rel0bar1} reads
 \be
  (\s_j\sb_j)\,\sb_i(\s_j\sb_j)^{-1}=\s_i\,,
 \ee
which follows from \eqref{rel0} in the form
$(\s_j\sb_j)\,\s_i(\s_j\sb_j)^{-1}=\sb_i$.

For $a_{ji}$ even we have to prove \eqref{rel0bar1}, which reads
 \be\label{rel0ev1}
  (\s_j\sb_j)\,\sb_i(\s_j\sb_j)^{-1}=\sb_i\,.
 \ee
This follows from \eqref{rel0},
$(\s_j\sb_j)\,\s_i(\s_j\sb_j)^{-1}=\s_i$, by multiplying both sides
by $\s_i\sb_i$ and using \eqref{etacomm}. $\Box$

\begin{prop}\label{EXT}
The group $W^U_G$ is given by the following extension
 \be\label{EXTU}
  1\longrightarrow \IZ_2^{|I|}\longrightarrow W^U_G\longrightarrow W_G
  \longrightarrow 1\,,
 \ee
of the Weyl group $W_G$ by the abelian group $\IZ_2^{|I|}\subset
W_G^U$.
\end{prop}

\proof Introduce the elements $\eta_i:=\s_i\sb_i,\,i\in I$.  They
have order two, and they pairwise commute by \eqref{etacomm}:
  \be\label{CONS1}
   \eta_i^2=1, \qquad
   \eta_i\eta_j=\eta_j\eta_i\,,\qquad i,j\in I\,,
   \ee
and are invariant under involution $\Pi$. Consider a $\Pi$-stable
subgroup $H_\eta\subset W^U_G$ generated by $\{\eta_i,\,i\in I\}$.
The relation \eqref{rel0} may be equivalently written in the
following form:
 \be\label{rel22}
  \s_i\eta_j\,=\,\eta_j\eta_i^{-a_{ji}}\s_i\,,\qquad
  \sb_i\eta_j\,=\,\eta_j\eta_i^{-a_{ji}}\sb_i\,.
  \ee
Indeed, from \eqref{rel0} we have
$\s_i\eta_j=\eta_j\Pi^{1-a_{ji}}(\sb_i)$, which for $a_{ji}$ even
reads $\s_i\eta_j=\eta_j \s_i $ and for $a_{ji}$ odd \eqref{rel0} is
equivalent to $\s_i\eta_j=\eta_j\sb_i=\eta_j\eta_i\s_i$. This
implies the first equation in \eqref{rel22}. The second relation in
\eqref{rel22} is obtained by applying the automorphism $\Pi$ to the
first one.
 The identities \eqref{rel22} yield that
the subgroup $H_\eta\subset W^U_G$ generated by $\eta_i,i\in I$ is
normal. A proof of the fact that $|H_\eta|=2^{|I|}$ is given in
Lemma \ref{ETA} below.

Next, we consider the quotient group $W_G^U/\IZ_2^{|I|}$. It is
generated by $s_i:=\pi(\s_i),\,i\in I$ and satisfying the standard
relations \eqref{BR0}, \eqref{BR} of the group $W_G$. Indeed,
$\pi(\s_i)=\pi(\bar{\s}_i)$ implies that the relations \eqref{rel}
are mapped to the relations \eqref{BR0}, relations \eqref{rel0}
equivalent to \eqref{rel22} become identities, and the braid
relations \eqref{Mrel} of $W_G^U$ are mapped to the braid relations
\eqref{BR} of $W_G$. Thus, we have a surjective homomorphism $\pi$
\be \pi\,:\quad
  W_G^U\,\longrightarrow\,W_G\,,\qquad\eta_i\,\longmapsto\,1\,,\quad
  i\in I\,.
  \ee
 This gives the exact sequence \eqref{EXTU}. $\Box$

The following analog of Theorem \ref{DTT} holds.

  \begin{te}\label{TE} Let $U\subset G$ be a maximal compact
    subgroup of the complex semisimple simply connected Lie group $G$. Let
$(\Delta,W_G)$ be the root system of
  $\Fg=\Lie(G)$ with the Cartan matrix $\|a_{ij}\|$. Let $\gamma$ be the
  generator of the   Galois group $\Gamma={\rm
    Gal}(\IC/\IR)$ of the field extension $\IR\subset\IC$ and let $\imath\in\IC$ be the imaginary unit.
Then the following map
 \be\label{GEN1}
  \sigma_i\longmapsto\varsigma_i:=e^{\frac{\imath\pi}{2}(e_i+f_i)} \,\gamma,\qquad
  \sb_i\longmapsto \bar{\varsigma}_i:=e^{-\frac{\imath\pi}{2}(e_i+f_i)}
  \,\gamma\,,\qquad i\in I,
 \ee
defines an injective homomorphism
$W^U_G\longrightarrow\widetilde{U}$, with $\widetilde{U}$ given by
\eqref{SDP}. The elements $\varsigma_i,\,i\in I$ and $\gamma$
together with the maximal torus $H$ generate the group
$N_G(H)\rtimes \Gamma$.
\end{te}
We give a proof of Theorem \ref{TE} in Section  \ref{PRR} below.

Let us stress a clear analogy between the constructions  of $W_G^T$
and $W_G^U$. On the one hand, in the Tits setting the finite group
is embedded into the maximal compact subgroup $K\subset G(\IR)$ of
the maximally split real form $G(\IR)\subset G$. On the other hand
the extension $W^U_G$ constructed above is embedded into Galois
group extension  $\tilde{U}$ of maximal compact subgroup $U\subset
G$ (note that the action of complex conjugation on $U$ may be
equivalently represented by the action of the Cartan anti-involution
\eqref{Cartaninv}). In this way $\tilde{U}$ plays the role analogous
to $K\subset G(\IR)$ in the Tits construction, while $W^U_G$ looks
like a ``complex'' analog of the finite group $W^T_G$, where the
relations $\tau_i^2=\theta_i$ are replaced by
$\sigma_i\bar{\sigma}_i=\eta_i$. Let us note that for each $i\in
  I$ we have the cyclic subgroup
  $\<\dot{s}_i:\,(\dot{s}_i)^4=1\>=\IZ_4\subset W_G^T$ in the Tits
  construction. Meanwhile in the case of $W_G^U$ for each $i\in I$ the
group $\<\s_i,\sb_i:\,\s_i^2=\sb_i^2=(\s_i\sb_i)^2=1\>=(\IZ_2)^2\subset
  W_G^U$ appears. Note that these two instances  exhaust possible
  extensions of $\IZ_2$ by   $\IZ_2$. It is natural to expect that with other real forms
of the complex group $G$ one can associate appropriate extensions of
the Weyl group $W_G$. These extensions presumably would be
combinations of both constructions considered above.

\section{Adjoint action of the extended Weyl groups}

While the action of $W_G$ on the maximal commutative subalgebra
$\mathfrak{h}={\rm Lie}(H)$ is defined canonically
\eqref{WeylGroupAction} and does not depend on a lift of $W_G$ into
$N_G(H)$ its action on the whole Lie algebra $\Fg={\rm Lie}(G)$ does
depend on the lift. Above we have considered two extensions of the Weyl group
$W_G$ together with their homomorphisms into the corresponding Lie
group. Here we describe their induced adjoint actions on $\Fg$.

\begin{prop}\label{PROP}
The adjoint action of the Tits group $W_G^T$ on the Lie algebra
$\Fg=\Lie(G)$ via  homomorphism \eqref{Tits} is given by
  \be\label{ADT1}
  \Ad_{\dot{s}_i}(e_i)\,=\,-f_i, \qquad
  \Ad_{\dot{s}_i}(f_i)\,=\,-e_i\, ,
  \ee
  \be \label{ADT11}
  \Ad_{\dot{s}_i}(e_j)=e_j,
  \qquad \Ad_{\dot{s}_i}(f_j)\,=\,f_j,\qquad a_{ij}=0\,,
  \ee
  \be \label{ADT2}
   \Ad_{\dot{s}_i}(e_j)\,
   =\,\frac{(-1)^{|a_{ij}|}}{|a_{ij}|!}
   \underbrace{\bigl[e_i,[\ldots[e_i}_{|a_{ij}|},\,e_j]\ldots]\bigr]\,,\\
   \Ad_{\dot{s}_i}(f_j)\,
   =\,\frac{1}{|a_{ij}|!}\underbrace{\bigl[f_i,[\ldots[f_i}_{|a_{ij}|},\,f_j]
   \ldots]\bigr],\qquad i\neq j\,.
    \ee
  \end{prop}

  \proof Relations \eqref{ADT1} are  actually  relations for
  $\mathfrak{sl}_2$ Lie subalgebras
generated by $\{e_i,h_i,f_i,\,i\in I\}$ and may easily be checked
using for example the standard faithful representation
\eqref{STREP}. Relations \eqref{ADT11} trivially follow from the Lie
algebra relations \eqref{Serre1}. Thus we need to prove
\eqref{ADT2}. Let us introduce the following notation:
$\dot{s}_i(a):=\Ad_{\dot{s}_i}(a)$. Then for the conjugated
generators we have
  \be
[h_k,\dot{s}_i(e_j)]=\dot{s}_i([h_{s_i(k)},e_j])=\<s_i(\alpha_k^{\vee}),\alpha_j\>\,
\dot{s}_i(e_j)=(a_{kj}-a_{ki}a_{ij})\dot{s}_i(e_j)\,,
\ee
\be
[h_k,\dot{s}_i(f_j)]=\dot{s}_i([h_{s_i(k)},f_j])=-\<s_i(\alpha_k^{\vee}),\alpha_j\>\,
\dot{s}_i(f_j)=-(a_{kj}-a_{ki}a_{ij})\dot{s}_i(f_j)\,.
\ee
 These relations fix the r.h.s. of \eqref{ADT2} up to coefficients. Let us calculate the
coefficients by taking into account  only the terms of the right
weights. We have
 \be
  \dot{s}_i(e_j)=e^{f_i}e^{-e_i}e^{f_i}\,e_j\,e^{-f_i}e^{e_i}e^{-f_i}=
  \frac{(-1)^{|a_{ij}|}}{|a_{ij}|!}\, \left({\rm
  Ad}_{e^{f_i}e_ie^{-f_i}}^{|a_{ij}|}(e_j)\right)+\cdots,
 \ee
where we have used the Serre relations \eqref{Serre2} and denote by
$\cdots$ the terms of the ``wrong'' weight.   Taking into account
 \be
  \Ad_{e^{f_i}}(e_i)=e_i+\cdots,
 \ee
we obtain the first relation in \eqref{ADT2}. The second relation is
obtained quite similarly  using the following equality (for a proof
see Lemma \ref{LEM0})
 \be\label{IDD}
 e^{f_i}e^{-e_i}e^{f_i}=e^{-e_i}e^{f_i}e^{-e_i}.
 \ee
In this case we have
 \be
  \dot{s}_i(f_j)=e^{-e_i}e^{f_i}e^{-e_i}\,f_j\,e^{e_i}e^{-f_i}e^{e_i}=
  \frac{1}{|a_{ij}|!}\,\left({\rm Ad}_{e^{-e_i}f_ie^{e_i}}^{|a_{ij}|}(f_j)\right)+ \cdots.
 \ee
Taking into account
 \be
  \Ad_{e^{-e_i}}(f_i)=f_i+\cdots,
 \ee
we obtain the second  relation in \eqref{ADT2}.  $\Box$

Let us stress that there is a simple way to get rid of sign factors  in
  \eqref{ADT1} and \eqref{ADT2}. Define a new set of generators
 $\tilde{e}_i=-e_i$, $\tilde{f}_i=f_i$. Then we have
 \be
  \Ad_{\dot{s}_i}(\tilde{e}_i)\,=\,\tilde{f}_i\,,
  \qquad \Ad_{\dot{s}_i}(\tilde{f}_i)\,=\,\tilde{e}_i\,,
  \ee
  \be
  \Ad_{\dot{s}_i}(\tilde{e}_j)=\tilde{e}_j\,,
  \qquad \Ad_{\dot{s}_i}(\tilde{f}_j)\,=\,\tilde{f}_j\,,\qquad a_{ij}=0\,,
  \ee
  \be
   \Ad_{\dot{s}_i}(\tilde{e}_j)\,
   =\,\frac{1}{|a_{ij}|!}
   \underbrace{\bigl[\tilde{e}_i,[\ldots[\tilde{e}_i}_{|a_{ij}|},\,\tilde{e}_j]
   \ldots]\bigr]\,,\\
   \Ad_{\dot{s}_i}(\tilde{f}_j)\,
   =\,\frac{1}{|a_{ij}|!}\underbrace{\bigl[\tilde{f}_i,[
     \ldots[\tilde{f}_i}_{|a_{ij}|},\,\tilde{f}_j]
   \ldots]\bigr],\qquad i\neq j\,.
    \ee

Now we describe the action on $\Fg$ of the Weyl group extension
$W_G^U$ introduced in Section 4. It is convenient to express it in
terms of purely imaginary generators $\imath e_i$, $\imath f_i$, $i\in
I$.

 \begin{prop} The elements of the group $W^U_G$ act on the Lie algebra
   $\Fg={\rm Lie}(G)$ via   homomorphism \eqref{GEN1} as follows
 \be\label{ADU2}
  \Ad_{\vs_i}(\imath e_i)\,=-\,\imath f_i\,,\qquad
  \Ad_{\vs_i}(\imath f_i)\,=-\,\imath e_i\,,
 \ee
  and
  \be \label{ADU21}
  \Ad_{\vs_i}(\imath e_j)=-\imath e_j,
  \qquad \Ad_{\vs_i}(\imath f_j)\,=\,-\imath f_j,\qquad a_{ij}=0\,,
  \ee
  \be\label{ADU1}
   \Ad_{\vs_i}(\imath e_j)\,
   =\,-\, \frac{1}{|a_{ij}|!}\,
   \underbrace{\bigl[\imath e_i,[\ldots[\imath e_i}_{|a_{ij}|},\,
   \imath e_j]\ldots]\bigr],\\
   \Ad_{\vs_i}(\imath f_j)\,
   =\,-\,\frac{1}{|a_{ij}|!}\,
   \underbrace{\bigl[\imath f_i,[\ldots[\imath f_i}_{|a_{ij}|},\,
   \imath f_j]\ldots]\bigr],\qquad i\neq j\,.
  \ee
\end{prop}

\proof Taking into account \eqref{Serre1} we have
 \be
  e^{\imath \pi t\, {\rm ad}_{h_i}}(e_j)=e_j\, e^{\imath \pi
  ta_{ij}}\, , \qquad
  e^{\imath \pi t\, {\rm ad}_{h_i}} (f_j) =f_j\, e^{-\imath \pi
  ta_{ij}}\, .
 \ee
Next, we use the following representation for $\vs_i$ (see Lemma
\ref{LEM0} for details):
 \be
  \vs_i=\dot{s}_i e^{\frac{\imath\pi}{2}h_i}\gamma\,,\qquad i\in I\,,
 \ee
and Proposition \ref{PROP} to obtain \eqref{ADU1} and \eqref{ADU2}.
$\Box$

\section{Proof of Theorem \ref{TE}} \label{PRR}

We start with establishing an explicit relation between the
generators $\vs_i,\vsb_i$ \eqref{GEN1} and the Tits generators
$\dot{s}_i$.

\begin{lem}\label{LEM0} For each $i\in I$ the following identities hold
 \be\label{ID1}
  \dot{s}_i:=e^{f_i}\,e^{-e_i}\,e^{f_i}=e^{-e_i}\,e^{f_i}\,e^{-e_i}
  = e^{\frac{\imath\pi}{4}h_i}\, e^{\frac{\imath\pi}{2}(e_i+f_i)}\,
  e^{-\frac{\imath\pi}{4}h_i},\\
  \dot{s}_i^2=e^{\imath \pi h_i}.
 \ee
Thus the generators $\{\vs_i,\vsb_i,\,i\in I\}$ defined by
\eqref{GEN1} may be represented as follows
 \be\label{ID2}
  \vs_i=e^{-\frac{\imath\pi}{2}h_i}\dot{s}_i\gamma
  =e^{-\frac{\imath\pi}{4}h_i}\, \dot{s}_i\,e^{\frac{\imath\pi}{4}h_i}\gamma,\qquad
  \vsb_i=e^{\frac{\imath\pi}{2}h_i}\dot{s}_i \gamma
  =e^{\frac{\imath\pi}{4}h_i}\, \dot{s}_i\,e^{-\frac{\imath\pi}{4}h_i}\gamma,\\
  \vs_i\vsb_i=\dot{s}_i^2=e^{\imath \pi h_i}\,.
 \ee
\end{lem}

\proof The identities  \eqref{ID1}  follow from the corresponding
relations in $SL_2\subset G$, using the standard faithful
two-dimensional representation $\phi: SL_2\to {\rm End}(\IC^2)$
 \be\label{STREP}
  \phi(e)=\begin{pmatrix} 0 &1 \\ 0& 0 \end{pmatrix}, \qquad
  \phi(f)=\begin{pmatrix} 0 &0 \\ 1&0 \end{pmatrix}, \qquad
  \phi(h)=\begin{pmatrix} 1 &0 \\ 0&-1 \end{pmatrix}.
 \ee
Direct calculations show that
 \be
  \phi(\dot{s})=\phi\left(e^{f}\,e^{-e}\,e^{f}\right)=\phi\left( e^{-e}\,e^{f}\,e^{-e}\right)=
  \begin{pmatrix} 0& -1 \\ 1 & 0  \end{pmatrix},
 \ee
 \be
  \phi(\dot{s}^2)
  =\phi\left( e^{-e}\,e^{f}\,e^{-e}\,e^{-e}\,e^{f}\,e^{-e}\right)=
  \begin{pmatrix} -1& 0 \\ 0 & -1  \end{pmatrix}
  =\phi\left(e^{\imath \pi h}\right),
 \ee
 \be
  \phi\left(e^{-\frac{\imath\pi}{4}h}\dot{s}e^{\frac{\imath\pi}{4}h}\right)=
  \phi\left(e^{\frac{\imath\pi}{2}(e+f)}\right)=
  \begin{pmatrix} 0& \imath \\ \imath & 0
    \end{pmatrix},\quad
  \phi\bigl(e^{\frac{-\imath\pi}{2}(e+f)}\bigr)=
  \begin{pmatrix} 0& -\imath \\-\imath & 0
    \end{pmatrix}.
 \ee
Then \eqref{ID1}, \eqref{ID2} follow from the faithfulness of
$\phi$. $\Box$

\begin{lem}   The following relations hold
 \be\label{SK}
  \vs_i^2=\vsb_i^2=1, \qquad i\in I\,.
 \ee
 \end{lem}
 \proof Direct calculation gives
 \be
  \vs_i^2= e^{-\frac{\imath\pi}{4}h_i} \dot{s}_i e^{\frac{\imath\pi}{4}h_i}\,\gamma
  e^{-\frac{\imath\pi}{4}h_i} \dot{s}_i e^{\frac{\imath\pi}{4}h_i}\,\gamma=
  e^{-\frac{\imath\pi}{4}h_i} \dot{s}_i e^{\frac{\imath\pi}{4}h_i}\,
  e^{\frac{\imath\pi}{4}h_i} \dot{s}_i e^{-\frac{\imath\pi}{4}h_i}\\=
  e^{-\frac{\imath\pi}{4}h_i} \dot{s}_i e^{\frac{\imath\pi}{2}h_i}\,
  \dot{s}_i e^{-\frac{\imath\pi}{4}h_i}
  = e^{-\frac{\imath\pi}{4}h_i}  e^{-\frac{\imath\pi}{2}h_i}\,e^{-\frac{\imath\pi}{4}h_i} \,
  (\dot{s}_i)^2=e^{-\imath \pi h_i}\cdot e^{\imath \pi h_i}=1\,.
 \ee
The identity $\vsb_i^2=1$ follows from $\vsb_i=\g\vs_i\g,i\in I$.
$\Box$

Now let us verify that the generators $\vs_i,\vsb_i$ and
$\xi_i=\vs_i\vsb_i,\,i\in I$ satisfy the remaining defining
relations \eqref{rel22}, \eqref{Mrel} for the group $W_G^U$:
 \be\label{rel2}
  \vs_i\xi_j\,=\,\xi_j\xi_i^{-a_{ji}}\vs_i\,,\qquad
  \vsb_i\xi_j\,=\,\xi_j\xi_i^{-a_{ji}}\vsb_i\,,
 \ee
and
 \be\label{braid}
  \underbrace{\vs_i\vs_j\vs_i\cdots}_{m_{ij}}\,=\,
  \underbrace{\vsb_j\vsb_i\vsb_j\cdots}_{m_{ij}}\,,
  \qquad i\neq j\,,\quad i,j\in I\,,
  \ee
where $m_{ij}=2,3,4,6$ for $a_{ij}a_{ji}=0,1,2,3$, respectively.

The first identity in \eqref{rel2} follows from \eqref{ID2} and
\eqref{TitsAd}:
 \be
  \vs_i\xi_j\vs_i^{-1}
  =e^{-\frac{\pi\imath}{2}h_i}\dot{s}_i\g
  e^{\pi\imath h_j}\g\dot{s}_i^{-1}e^{\frac{\pi\imath}{2}h_i}
  =e^{-\frac{\pi\imath}{2}h_i}e^{\pi\imath(h_j-a_{ji}h_i)}
  e^{\frac{\pi\imath}{2}h_i}=\xi_j\xi_i^{-a_{ji}}.
 \ee
The other identity in \eqref{rel2} follows from
$\vsb_i=\g\vs_i\g,i\in I$.

For the relation \eqref{braid}, on the left side we have
 \be
  \underbrace{\vs_i\vs_j\vs_i\cdots}_{m_{ij}}
  =\bigl(\underbrace{e^{-\frac{\pi\imath}{2}h_i}\dot{s}_i
  e^{\frac{\pi\imath}{2}h_j}\dot{s}_j\cdots}_{m_{ij}}\bigr)\g^{m_{ij}}\\
  =\exp\frac{\pi\imath}{2}\bigl(\underbrace{-h_i+s_ih_j-s_is_jh_i+\ldots}_{m_{ij}}\bigr)
  (\underbrace{\dot{s}_i\dot{s}_j\cdots}_{m_{ij}})\g^{m_{ij}},
 \ee
and on the right side:
 \be
  \underbrace{\vsb_j\vsb_i\vsb_j\cdots}_{m_{ij}}
  =\bigl(\underbrace{e^{\frac{\pi\imath}{2}h_j}\dot{s}_j
  e^{-\frac{\pi\imath}{2}h_i}\dot{s}_i\cdots}_{m_{ij}}\bigr)\g^{m_{ij}}\\
  =\exp\frac{\pi\imath}{2}\bigl(\underbrace{h_j-s_jh_i+s_js_ih_j-\ldots}_{m_{ij}}\bigr)
  (\underbrace{\dot{s}_j\dot{s}_i\cdots}_{m_{ij}})\g^{m_{ij}}.
 \ee
Since
$\underbrace{\dot{s}_i\dot{s}_j\cdots}_{m_{ij}}=\underbrace{\dot{s}_j\dot{s}_i\cdots}_{m_{ij}}$
holds due to \eqref{Tits1}, the identity \eqref{braid} reduces to
the following
 \be\label{EQQ}
  \exp\frac{\pi\imath}{2}\bigl(\underbrace{-h_i+s_ih_j-s_is_jh_i+\ldots}_{m_{ij}}\bigr)=
  \exp\frac{\pi\imath}{2}\bigl(\underbrace{h_j-s_jh_i+s_js_ih_j-\ldots}_{m_{ij}}\bigr).
 \ee
In turn the identity \eqref{EQQ} may be proved by invoking the
following fact.

\begin{lem}\label{Dihedral}
For each pair $i,j\in I$, $i\neq j$, the following holds:
 \be\label{Dm}
  (\underbrace{1-s_j+s_is_j-\ldots}_{m_{ij}})h_i=0\,.
 \ee
\end{lem}
\proof Consider the order $2m_{ij}$ Coxeter subgroup of $W_G$
generated by a pair of the simple root reflections $s_i,s_j,\,i\neq
j$\,:
 \be
  \<s_i,s_j:\,(s_is_j)^{m_{ij}}=s_i^2=s_j^2=1\>\subset W_G.
 \ee
This group is isomorphic to the dihedral group $D_{m_{ij}}\subset
O_2(\IR)$ of symmetries of $m_{ij}$-gone in the real plane
$V_{ij}=\IR h_i\oplus\IR h_j$. The dihedral group may be
equivalently written in the following form:
 \be
  D_{m_{ij}}=\<t,r:\,t^{m_{ij}}=r^2=1,\,\,
  rtr^{-1}=t^{-1}\>,\quad t=s_is_j,\,\,r=s_i\,.
 \ee
We have two projectors in the plane $V_{ij}=\IR h_i\oplus\IR h_j$:
 \be
  P_{\pm}=\frac{1\pm s_i}{2}:\,V_{ij}\longrightarrow V_{ij},\qquad
  P_{\pm}^2=P_{\pm}, \quad P_{\pm} P_{\mp}=0,
 \ee
such that
 \be
  P_+h_i=0, \qquad P_-h_i=h_i\,.
 \ee
Therefore, the identity \eqref{Dm} is equivalent to the following:
 \be\label{DmId}
  (\underbrace{1-s_j+s_is_j-\ldots}_{m_{ij}})(1-s_i)\,\cdot h_i
  =\sum_{g\in D_{m_{ij}}}(-1)^gg\cdot h_i=0,
 \ee
where $(-1)^g:=\det(g)$ is the sign character of $D_{m_{ij}}\subset
O_2(\IR)$. The kernel of the sign character is a normal subgroup,
 \be
  \IZ_{m_{ij}}=\<t:\,t^{m_{ij}}=1\>\subset D_{m_{ij}},
 \ee
which consists of the rotations by $\frac{2\pi k}{m_{ij}},\,0\leq
k<m_{ij}$ of the plane $V_{ij}$. The non-trivial co-set in
$D_{m_{ij}}/\<t\>$ consists of the reflections $\{r_k=t^kr:\,0\leq
k<m_{ij}\}$ with $r_0=r$ being a reflection sending $h_i$ to $-h_i$.
Thus we have $(-1)^{r_k}=\det(r_k)=-1$ and
 \be
  D_{m_{ij}}=\{t^k:\,0\leq k<m_{ij}\}\sqcup\{r_k=rt^k:\,0\leq
  k<m_{ij}\},
 \ee
hence the identity \eqref{DmId} reads
 \be
  \sum_{g\in D_{m_{ij}}}(-1)^gg\cdot h_i
  =(1-r)\sum_{k=0}^{m_{ij}-1} t^k\cdot h_i.
 \ee
Now in the group algebra $\IC[D_{m_{ij}}]$ the following identity
holds:
 \be
  t^{m_{ij}}-1=(t-1)\sum_{k=0}^{m_{ij}-1} t^k=0\,.
 \ee
Since $t$ acts in the faithful representation $V_{ij}$ without fixed
vectors, we infer that $\sum\limits_{k=0}^{m_{ij}-1}t^k \cdot h_i=0$
and thus prove \eqref{DmId}.\, $\Box$

\begin{lem}\label{ETA} The elements $\eta_i=\s_i\sb_i,\,i\in I$
generate a subgroup $H_\eta\subset W_G^U$ of order
  $|H_\eta|=2^{|I|}$.
\end{lem}
\proof By \eqref{SK}, \eqref{rel2} and \eqref{braid}, the elements
$\{\vs_i,\vsb_i,\,i\in I\}\subset\tilde{U}=U\rtimes\Gamma$ satisfy
the defining relations of the group $W_G^U$ from Definition
\ref{DDD}. Moreover, the images $\vs_i\vsb_i=e^{\imath\pi h_i}\in U$
of the elements $\eta_i=\s_i\sb_i,\,i\in I$ generate the subgroup
$H_{\eta}\simeq H^{(2)}$ of order two points in the maximal torus
$H\subset G$, so that the order of $H_\eta$ should be not less then
$2^{|I|}$ and thus $|H_\eta|=2^{|I|}$ holds. $\Box$

We complete our proof of Theorem \ref{TE} by verifying injectivity
of the homomorphism  $\psi:\,W^{U}_G\to \tilde{U}$. By Proposition
\eqref{EXT}, $W^{U}_G$ has a structure of the group extension:
 \be\label{WGUext}
  1\longrightarrow
  \IZ_2^{|I|}\longrightarrow W_G^U\longrightarrow W_G\longrightarrow1.
 \ee
Let $\CW_G^U=\psi(W_G^U)\subset\tilde{U}$, then the $\CW_G^U$-action
on $\mathfrak{h}=\Lie(H)$ implies the existence of the surjective
homomorphism $\pi:\,\CW_G^U\to W_G$. By Lemma \ref{ETA}, $\CW_G^U$
contains a normal abelian subgroup generated by
$\{\psi(\eta_i)=\vs_i\vsb_i=e^{\imath\pi h_i},\,i\in I\}\subset U$,
which is isomorphic to $H^{(2)}\simeq\IZ_2^{|I|}$. Clearly, the
normal abelian subgroup $H^{(2)}$ acts trivially on $\fh$, hence it
is in the kernel of the surjective homomorphism $\pi$, which entails
$|W_G^U|\leq |\CW_G^U|$. On the other hand, the existence of
homomorphism $\psi:\, W_G^U \to \CW_G^U$ implies that $|W_G^U|\geq
|\CW_G^U|$ and hence $|W_G^U|=|\CW_G^U|$. Thus for $\CW_G^U$ we have
the following exact sequence
 \be
  1\longrightarrow
  \IZ_2^{|I|}\longrightarrow \CW_G^U\longrightarrow
  W_G\longrightarrow1\,.
 \ee
Taking into account \eqref{WGUext} this provides a proof of
injectivity of $\psi$, and therefore, of Theorem \ref{TE}.

\noindent {\small {\bf A.A.G.} {\sl Laboratory for Quantum Field Theory
and Information},\\
\hphantom{xxxx} {\sl Institute for Information
Transmission Problems, RAS, 127994, Moscow, Russia};\\
\hphantom{xxxx} {\sl Interdisciplinary Scientific Center
  J.-V. Poncelet (CNRS UMI 2615),\\
 \hphantom{xxxx}  Independent University of Moscow,  Moscow, Russia;}\\
\hphantom{xxxx} {\it E-mail address}: {\tt anton.a.gerasimov@gmail.com}}\\
\noindent{\small {\bf D.R.L.}
{\sl Laboratory for Quantum Field Theory
and Information},\\
\hphantom{xxxx}  {\sl Institute for Information
Transmission Problems, RAS, 127994, Moscow, Russia};\\
\hphantom{xxxx} {\sl Moscow Center for Continuous Mathematical
Education,\\
\hphantom{xxxx} 119002,  Bol. Vlasyevsky per. 11, Moscow, Russia};\\
\hphantom{xxxx} {\it E-mail address}: {\tt lebedev.dm@gmail.com}}\\
\noindent{\small {\bf S.V.O.} {\sl
 School of Mathematical Sciences, University of Nottingham\,,\\
\hphantom{xxxx} University Park, NG7\, 2RD, Nottingham, United Kingdom};\\
\hphantom{xxxx} {\sl
 Institute for Theoretical and Experimental Physics,
117259, Moscow, Russia};\\
\hphantom{xxxx} {\it E-mail address}: {\tt oblezin@gmail.com}}

\end{document}